\documentclass{gtart}

\def\ifplaintex{\expandafter\ifx\csname documentclass\endcsname\relax}

\def\gtp{{\mathsurround=0pt\it $\cal G\mskip-2mu$eometry \&\ 
$\cal T\!\!$opology $\cal P\!$ublications}}  

\def\recd{{\small Received:\qua\receiveddate\ifx\reviseddate\relax
\else\qquad Revised:\qua\reviseddate\fi\par}} 

\def\lognumber#1{\def\thelognumber{#1}}
\def\volumenumber#1{\def\thevolumenumber{#1}}
\def\volumeyear#1{\def\thevolumeyear{#1}}
\def\papernumber#1{\def\thepapernumber{#1}}
\def\pagenumbers#1#2{\def\startpage{#1}\def\finishpage{#2}}
\def\published#1{\def\publishdate{#1}}

\def\received#1{\def\receiveddate{#1}}
\def\revised#1{\def\reviseddate{#1}}
\def\accepted#1{\def\accepteddate{#1}}
\def\asciititle#1{\def\theasciititle{#1}}

\long\def\asciiabstract#1{\long\def\theasciiabstract{#1}}
\def\asciikeywords#1{\def\theasciikeywords{#1}}

\let\\\par\let\thelognumber\relax\let\thevolumenumber\relax
\let\thepapernumber\relax\let\thevolumeyear\relax\let\startpage\relax
\let\finishpage\relax\let\publishdate\relax\let\receiveddate\relax
\let\reviseddate\relax\let\accepteddate\relax\let\theasciititle\relax
\let\theasciiauthors\relax
\let\theasciiabstract\relax\let\theasciikeywords\relax

\let\theasciiemail\relax

\ifplaintex
\font\logobig=cmssbx10 scaled 3836
\font\logomed=cmssbx10 scaled 2557
\else
\font\logobig=cmssbx10 scaled 4200
\font\logomed=cmssbx10 scaled 2800
\fi

\long\def\makeagttitle{   
\count0=\startpage
\agt\hfill      
\hbox to 45truept{\vbox to 0pt{\vglue -13truept{\logomed A\kern -.37em{\logobig 
T}\kern -.38em G}\vss}\hss}
\break
{\small Volume \thevolumenumber\ (\thevolumeyear)
\startpage--\finishpage\nl
Published: \publishdate}

\vglue .25truein

{\parskip=0pt\leftskip 0pt plus
1fil\def\\{\par\smallskip}{\Large\bf\thetitle}\par\medskip} \vglue
0.05truein

{\parskip=0pt\leftskip 0pt plus 1fil\def\\{\par}{\sc\theauthors}
\par\medskip}%
 
\vglue 0.03truein 

{\small\leftskip 25truept\rightskip 25truept{\bf Abstract}\stdspace\theabstract

{\bf AMS Classification}\stdspace\theprimaryclass
\ifx\thesecondaryclass\relax\else; \thesecondaryclass\fi\par
{\bf Keywords}\stdspace \thekeywords\par}\vglue 7truept

}   

\ifplaintex
\font\phead=cmsl9 scaled 950
\font\pnum=cmbx10 scaled 913
\font\pfoot=cmsl9 scaled 950
\headline{\vbox to 0pt{\vskip -4.5mm\line{\small\phead\ifnum
\count0=\startpage ISSN 1472-2739 (on-line) 1472-2747 (printed)
\hfill {\pnum\folio}\else\ifodd\count0\def\\{ }%
\ifx\theshorttitle\relax\thetitle\else\theshorttitle\fi\hfill{\pnum\folio}
\else\def\\{ and }{\pnum\folio}\hfill\ifx\theshortauthors\relax\theauthors
\else\theshortauthors\fi\fi\fi}\vss}}
\footline{\vbox to 0pt{\vglue 0mm\line{\small\pfoot\ifnum\count0=\startpage
\copyright\ \gtp\hfill\else
\agt, Volume \thevolumenumber\ (\thevolumeyear)\hfill\fi}\vss}}
\else
\font=cmsl9 scaled 1050
\font\lnum=cmbx10 
\font=cmsl9 scaled 1050
\makeatletter
\def\@oddhead{{\small\lhead\ifnum\count0=\startpage ISSN 1472-2739 
(on-line) 1472-2747 (printed)\hfill {\lnum\number\count0}\else\ifodd\count0
\def\\{ }\ifx\theshorttitle\relax \thetitle \else\theshorttitle\fi\hfill
{\lnum\number\count0}\else\def\\{ and }{\lnum\number\count0}
\hfill\ifx\theshortauthors\relax 
\theauthors\else\theshortauthors\fi\fi\fi}}\def\@evenhead{\@oddhead}
\def\@oddfoot{\small\lfoot\ifnum\count0=\startpage\copyright\ \gtp\hfill\else
\agt, Volume \thevolumenumber\ (\thevolumeyear)\hfill\fi}
\def\@evenfoot{\@oddfoot}
\makeatother
\fi
\let\maketitlepage\makeagttitle
\let\makeshorttitle\maketitlepage
\let\maketitle\maketitlepage

\newwrite\gtoutfile
\long\gdef\makeheadfile{  
{\def\\{, }\def\s{ }
\immediate\openout\gtoutfile head.xxx
\immediate\write\gtoutfile{Proxy-for: \ifx\theasciiauthors\relax
\theauthors\else\theasciiauthors\fi\s<\ifx\theasciiemail\relax\theemail\else\theasciiemail\fi>}
\immediate\write\gtoutfile{\noexpand\\}
\immediate\write\gtoutfile{Authors: \ifx\theasciiauthors\relax
\theauthors\else\theasciiauthors\fi}
{\def\\{ }\immediate\write\gtoutfile{Title: \ifx\theasciititle\relax
\thetitle\else\theasciititle\fi}}
\immediate\write\gtoutfile{Subj-class: GT or SG, GR etc}
\immediate\write\gtoutfile{MSC-class: \theprimaryclass\ifx\thesecondaryclass\relax\else, \thesecondaryclass\fi}
\immediate\write\gtoutfile{Journal-ref: Algebr. Geom. Topol. \thevolumenumber\s
(\thevolumeyear) \startpage-\finishpage}
\immediate\write\gtoutfile{Comments: Published by Algebraic and
Geometric Topology at}
\immediate\write\gtoutfile{\s\s\s  http://www.maths.warwick.ac.uk/agt/AGTVol\thevolumenumber/agt-\thevolumenumber-\thepapernumber.abs.html}
\immediate\write\gtoutfile{\noexpand\\}
\immediate\write\gtoutfile{}
\ifx\theasciiabstract\relax
\immediate\write\gtoutfile{\theabstract}\else
\immediate\write\gtoutfile{\theasciiabstract}\fi
\immediate\write\gtoutfile{}
\immediate\write\gtoutfile{\noexpand\\}
\immediate\write\gtoutfile{}
\immediate\closeout\gtoutfile}}  

\def\maketitlepage{\makeagttitle\makeheadfile}
\let\makeshorttitle\maketitlepage
\let\maketitle\maketitlepage

\lognumber{2}
\volumenumber{4}
\volumeyear{2004}
\papernumber{2}
\published{17 January 2004}
\pagenumbers{23}{29}
\received{2 September 2003}
\revised{3 October 2003}
\accepted{19 December 2003}

\usepackage{amsmath,amssymb,graphicx,subfigure}

\newtheorem{thm}{Theorem}[section]    
\theoremstyle{definition}
\newtheorem{defn}[thm]{Definition}   

\newcommand{\Cannon}{MR88a:20049}
\newcommand{\Epstein}{MR93i:20036}
\newcommand{\IndiraKim}{math.GR/0310356}
\newcommand{\Thiel}{MR95e:20052}
\newcommand{\Bowditch}{MR96b:20046}
\newcommand{\Olshan}{MR93d:20067}
\newcommand{\BridsonHaefliger}{MR1744486}
\newcommand{\BB}{MR2001j:20046}
\newcommand{\Gerstennotes}{MR2000a:20093}
\newcommand{\NSnotes}{NSnotes}

\newcommand{\bi}{\begin{itemize}}
\newcommand{\ei}{\end{itemize}}
\newcommand{\be}{\begin{enumerate}}
\newcommand{\ee}{\end{enumerate}}
\newcommand{\bc}{\begin{center}}
\newcommand{\ec}{\end{center}}
\newcommand{\bt}{\begin{tabular}}
\newcommand{\et}{\end{tabular}}
\newcommand{\smallcaps}[1]{\textrm{\textsc{#1}}}

\newcommand{\cat}{\smallcaps{CAT}}

\newcommand{\ac}{almost convex}
\newcommand{\gset}{generating set}
\newcommand{\iif}{isoperimetric function}

\newcommand{\cg}{Cayley graph}

\newcommand{\Z}{\mathbb Z}

\begin{document}

\title{$L_\delta$ groups are almost convex and have a\\sub-cubic Dehn function}
\asciititle{L_delta groups are almost convex and have a sub-cubic 
Dehn function}                    
\authors{Murray Elder}                  
\address{School of Mathematics and Statistics\\
University of St.\ Andrews\\
North Haugh, St.\ Andrews\\
Fife, KY16 9SS, Scotland}                  
\email{murray@mcs.st-and.ac.uk}                     
\url{http://turnbull.mcs.st-and.ac.uk/\char'176murray/}

\begin{abstract}
We prove that if the \cg\ of a finitely generated group enjoys the
property $L_\delta$  then
the group is \ac\ and  has a sub-cubic isoperimetric function.
\end{abstract}

\asciiabstract{%
We prove that if the Cayley graph of a finitely generated group enjoys
the property L_delta then the group is almost convex and has a
sub-cubic isoperimetric function.}

\primaryclass{20F65}                
\secondaryclass{20F67}                                
\keywords{Almost convex, isoperimetric function, property $L_\delta$}
\asciikeywords{Almost convex, isoperimetric function, property L_delta}

\maketitle

\section{Definitions}
In this article we show how a new metric property of groups called
$L_{\delta}$ is related to some older metric properties of groups.  We
prove that if a group has a finite generating set that enjoys $L_{\delta}$
then the group is almost convex with respect to this generating set, and
has a sub-cubic isoperimetric (or Dehn) function.  My thanks to Kim Ruane
and Indira Chatterji for introducing me to $L_{\delta}$, to Andrew Rechnitzer
for help with the figures, and to an anonymous reviewer for helpful
suggestions.

Let $(X,d)$ be a weakly geodesic metric space (in this paper take the
\cg\ with the word-metric).

\begin{defn}[$\delta$-path]
For any $\delta \geq 0$ and finite sequence of points $x_1,\ldots ,
x_n$, we say $(x_1,\ldots ,x_n)$ is a {\em $\delta$-path} if
\begin{eqnarray*}
d(x_1,x_2)+\ldots + d(x_{n-1},x_n)\leq d(x_1,x_n)+\delta.
\end{eqnarray*}
\end{defn}

\begin{defn}[$L_\delta$]
For any $\delta \geq 0$ the space $X$ has {\em property $L_\delta$} if
for each three distinct points $x,y,z\in X$ there exists a point $t\in
X$ so that the paths $(x,t,y),(y,t,z)$ and $(z,t,x)$ are all $\delta$-paths.
\end{defn}

\begin{figure}
\bt{ccc}
\bt{c}\includegraphics[width=5.5cm]{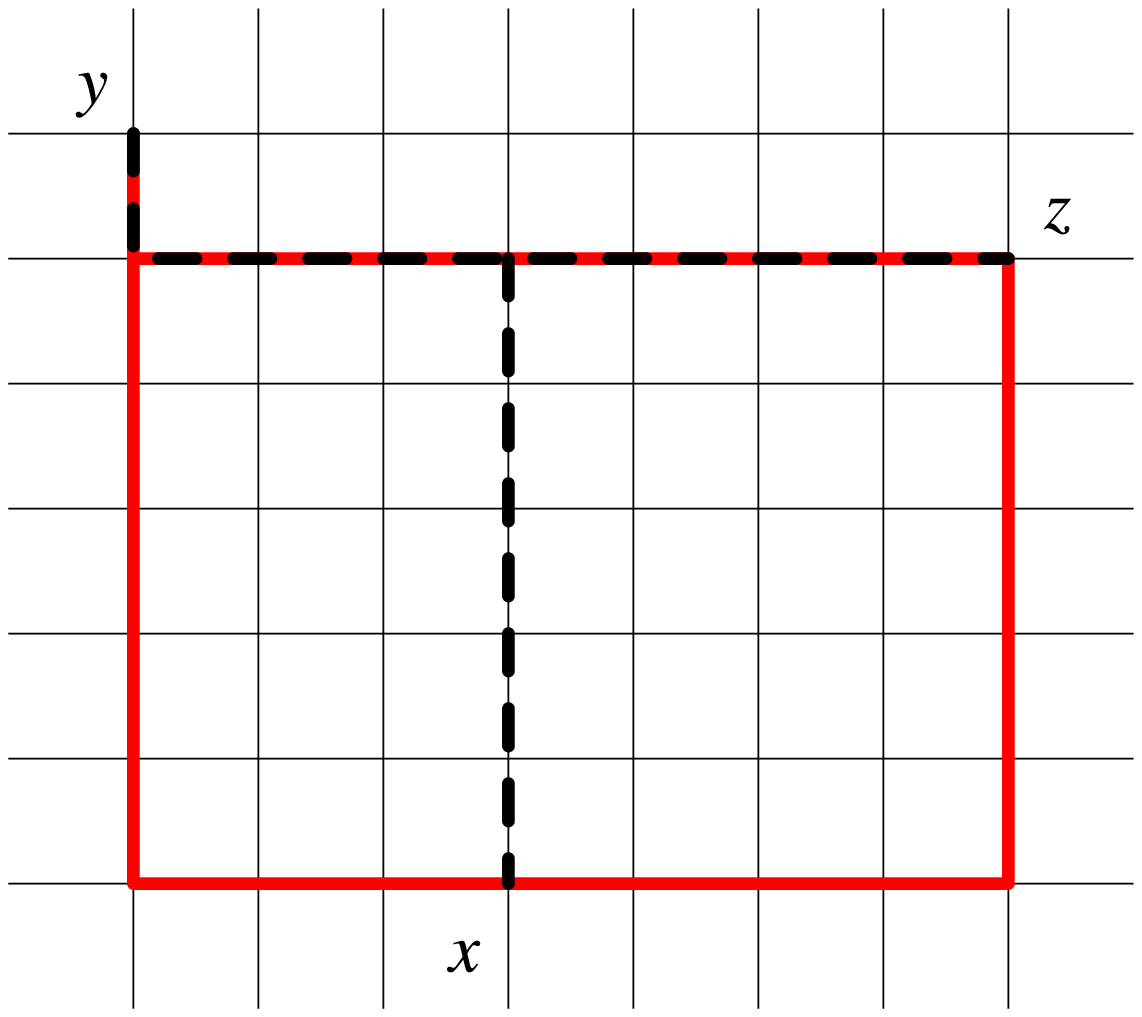}\et& &
\bt{c}\includegraphics[width=5.5cm]{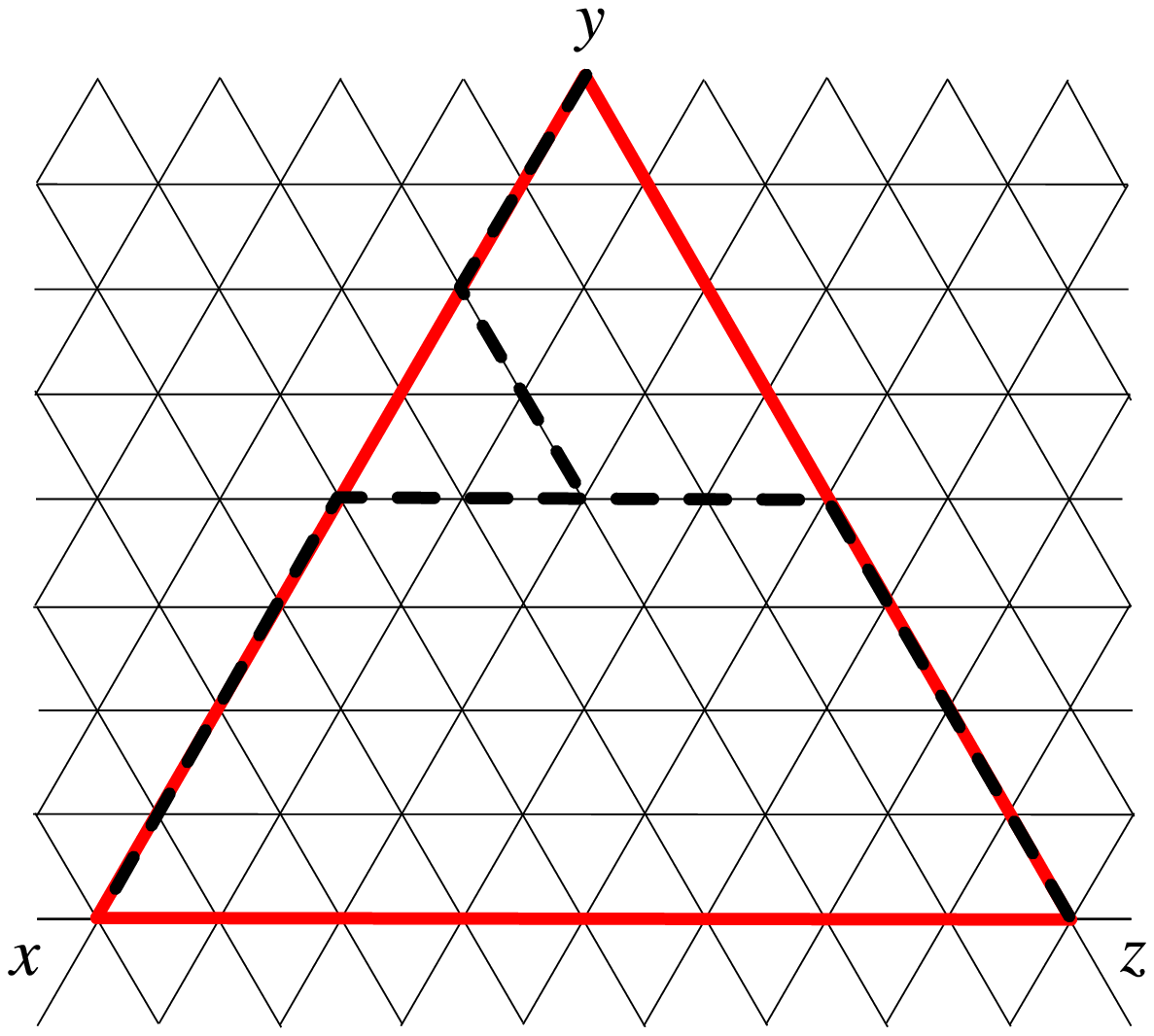}\et
\et
\caption{$\Z^2$ is $L_{\delta}$ for one \gset\ and not another}
\label{fig:abcPres}
\end{figure}

Examples of groups enjoying this property are word-hyperbolic groups,
 fundamental groups of various cusped hyperbolic $3$-manifolds,
groups that act on\break\eject 
\cat(0) cube complexes and products of trees, and
 Coxeter groups \cite{\IndiraKim}. The property is related to the
 property of {\em rapid decay} and the Baum-Connes conjecture. The property is not invariant under change of
 finite generating set. For example, Chatterji and Ruane show that
 $\Z^2$ with the usual generating set has $L_{\delta}$, however the
 generating set $\langle a,b,c \; | \; ab=ba=c\rangle$ does not enjoy
 the property, as illustrated in Figure \ref{fig:abcPres}.

Define the closed metric ball of radius $n$ to be the set of points
that lie within distance $n$ of the identity vertex in the \cg.
\begin{defn}[Almost convex]
A group $G$ with finite generating set $\mathcal G$ is {\em \ac} if
there is a constant $C\geq 0$ so that for any two vertices that lie
distance $n\geq 0$ from the identity vertex and at most 2 apart from
each other, there is a path connecting them that lies inside the
closed ball of radius $n$ and has length at most $C$.
\end{defn}
Cannon \cite{\Cannon} showed that if a group has this property then
one has an algorithm to construct any finite portion of the \cg. It
also implies finite presentability. The property is dependent on
choice of finite generating set \cite{\Thiel}.

Suppose that $G$ is a finitely presented group, with inverse closed
finite generating set $\mathcal G$ and finite set of relators
$\mathcal R\subseteq \mathcal G^*$.  Let $F(\mathcal G)$ be the free
group generated by $\mathcal G$.
 A word in $\mathcal G^*$ represents the identity in $G$ if and only
if it is freely equal to an expression of the form
$$\Pi_{i=1}^k g_ir_ig_i^{-1}$$ where the $g_i\in F(\mathcal G)$ and
$r_i\in \mathcal R\cup \mathcal R^{-1}$.

\begin{defn}[van Kampen diagram, Area]
Let $\Delta$ be a labeled, simply connected planar $2$-complex,
such that each edge is oriented and labeled by and element of
$\mathcal G$, and reading the labels on the boundary of each $2$-cell
gives and element of $\mathcal R\cup \mathcal R^{-1}$. We say 
$\Delta$ is a
{\em van Kampen diagram} for $w$ if reading the labels around the
boundary of $\Delta$ gives $w$.

Each van Kampen diagram with $k$ $2$-cells for $w$ gives a way of
expressing $w$ as a product $\Pi_{i=1}^k g_ir_ig_i^{-1}$. Conversely, each product
$\Pi_{i=1}^k g_ir_ig_i^{-1}$ gives a van Kampen diagram for $w$ with
at most $k$ $2$-cells. Thus van Kampen's Lemma states that $w$ has a van Kampen diagram if and only
if $w$ represents the identity element.  Define the {\em area} $A(w)$
of a word $w\in \mathcal G^*$ which represents the identity to be the
minimum $k$ in any such expression for $w$, or equivalently the
minimum number of $2$-cells in a van Kampen diagram for $w$.
\end{defn}

 In terms of the \cg, if a word evaluates to the identity element it
corresponds to a closed path, and its area is the least number of
relators needed to fill in this closed path.
See \cite{\BridsonHaefliger} p.155  for more details on van Kampen's Lemma.

\begin{defn}[Dehn function, Isoperimetric function]
 The {\em Dehn function} for $\langle \mathcal G|\mathcal R\rangle$ is
defined to be $\mathcal D(n)=\max\{A(w):w$ has at most $n$ letters,
and $w$ evaluates to the identity$\}$.  An {\em isoperimetric function}
for $\langle \mathcal G |\mathcal R\rangle$ is any function which
satisfies $f(n)\geq \mathcal D(n)$.
\end{defn}

\noindent
Two functions $f,g$ are said to be equivalent if there are constants
$A,A',B,B',$ $C,C',D, D',E,E'$ so that $f(n)\leq Ag(Bn+C)+Dn+E$, and $g(n)\leq
A'f(B'n+C')+D'n+E'$.  Up to this notion of equivalence, a Dehn
function for a group is invariant of the finite presentation.  If $G$
has a sub-quadratic \iif\ then it has a linear \iif\
\cite{\Bowditch,\Olshan}. For more details about isoperimetric
functions and van Kampen diagrams, see for example
\cite{\NSnotes,\Gerstennotes}.


\section{Results}

\begin{thm}\label{thm:ac}
If the \cg\ for a group has property $L_\delta$ for some $\delta \geq
0$ then it is \ac\ with constant $3\delta+2$.
\end{thm}

\begin{proof}
Consider two vertices $g,g'$ that lie at distance $n$ from the
identity in the \cg, and distance at most 2 apart.  We will call the
identity vertex $z$. Let $w$ and $w'$ be geodesic paths of length $n$
from $z$ to $g$ and $g'$ respectively.

If $n\leq \delta$ then the path $w^{-1}w'$ lies inside the closed ball
of radius $n$, connects $g$ to $g'$ and has length $2n\leq 2\delta$.

If $n>\delta$ then let $x$ be the point that lies distance exactly
  $n-\frac{\delta}2$ from $z$ along $w$, and $y$ the point that lies
  distance $n-\frac{\delta}2$ from $z$ along $w'$. Note that $x$ and
  $y$ need not be vertices (they may lie in the interior of
  edges). See Figure~\ref{fig:ac_1}.

\begin{figure}
\bt{ccc}
\bt{c}
\includegraphics[width=5.5cm]{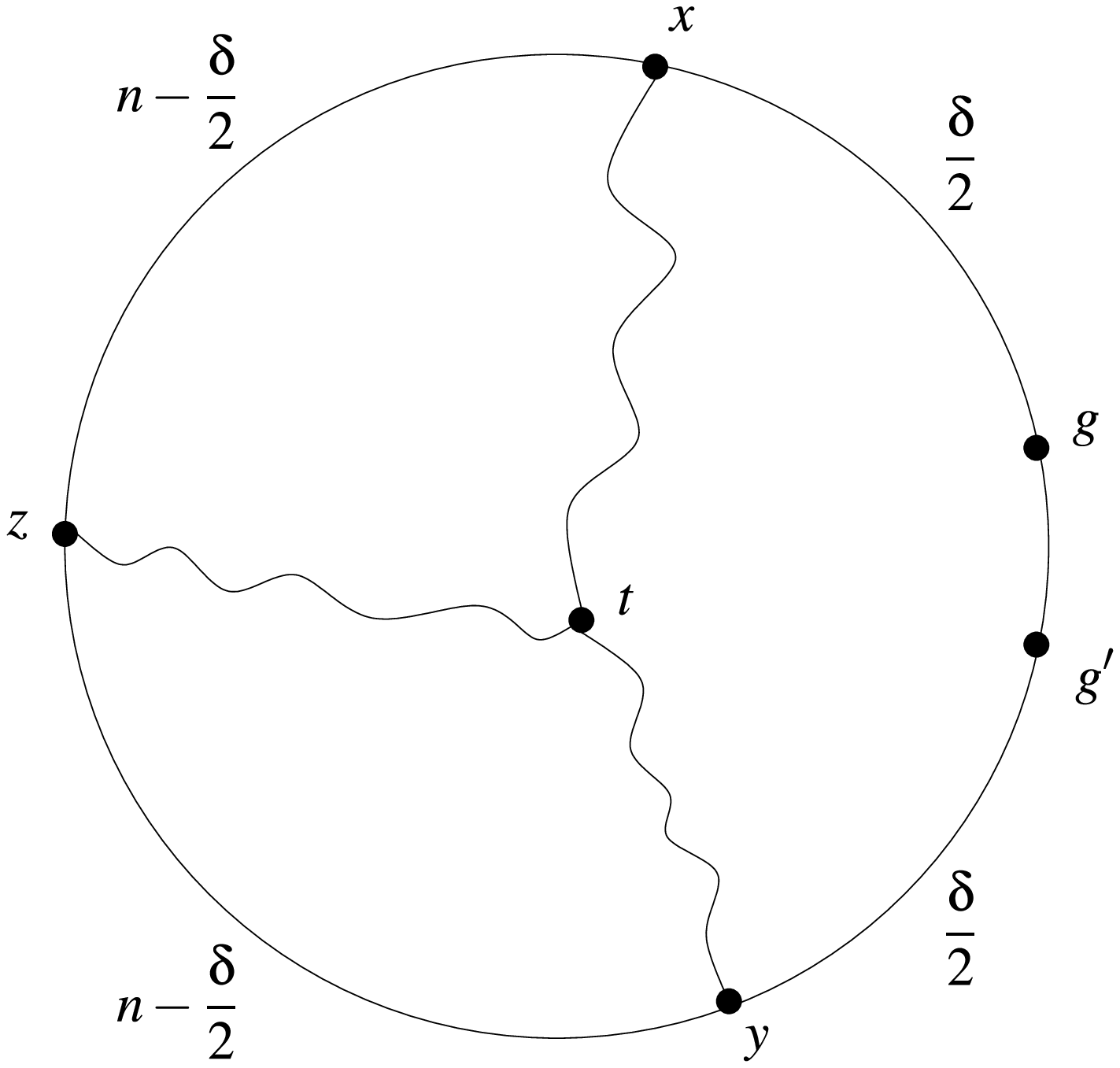}
\et& &
\bt{c}
\includegraphics[width=5.5cm]{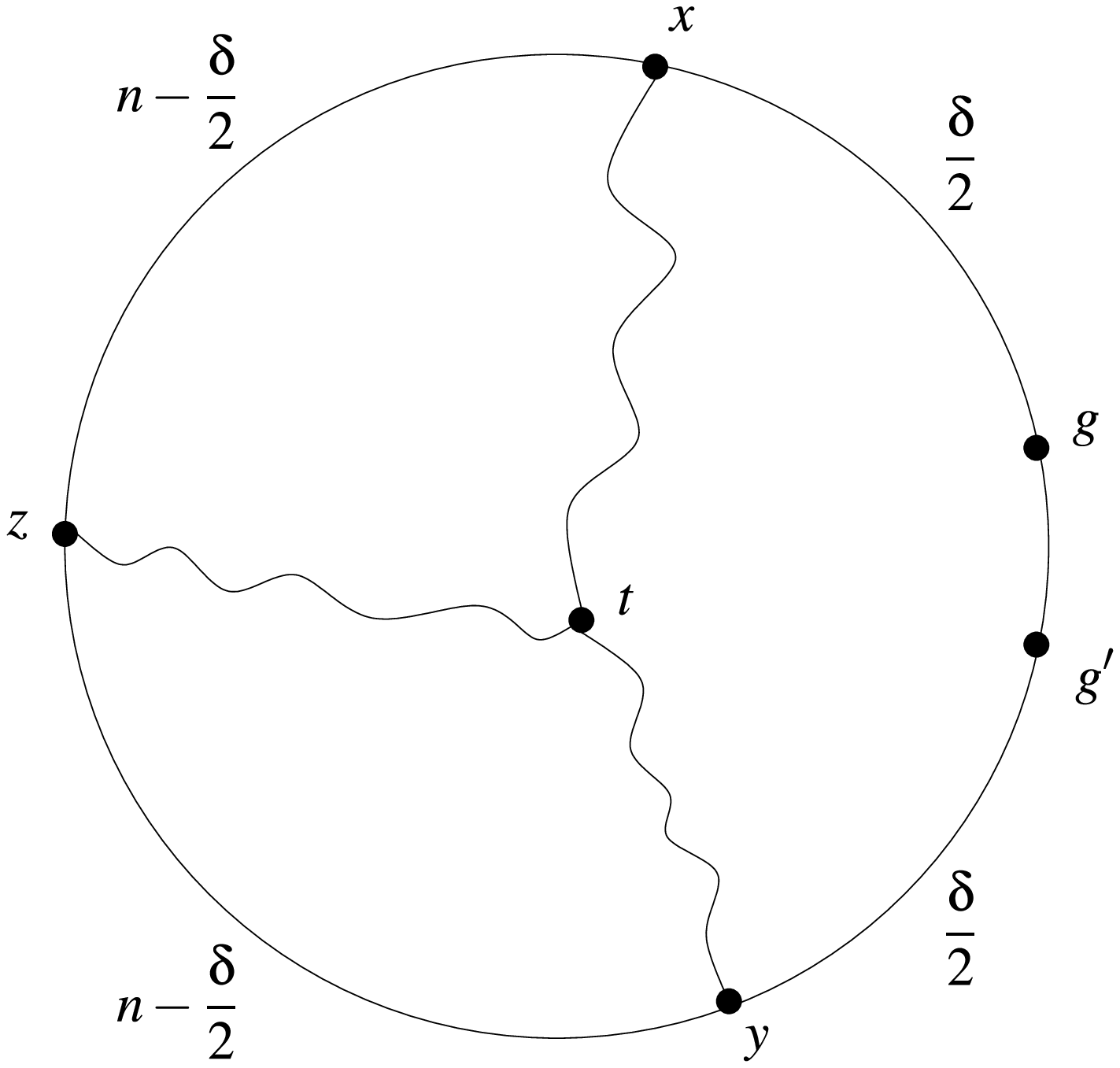}
\et
\et
\caption{$L_{\delta}$ implies \ac}
\label{fig:ac_1}
\end{figure}

By property $L_{\delta}$ there exists a point $t$ such that $d(z,t)+d(t,x)\leq
d(z,x)+\delta = n-\frac{\delta}2 + \delta =n+\frac{\delta}2$ and
$d(z,t)+d(t,y)\leq d(z,y)+\delta = n-\frac{\delta}2 + \delta
=n+\frac{\delta}2$ since $x$ and $y$ lie on geodesics. It follows
that the geodesic paths from $x$ to $t$ and $t$ to $y$ lie in the
closed ball of radius $n$.

Also, $d(x,t)+d(t,y)\leq d(x,y)+\delta \leq (\frac{\delta}2 + 2 +
\frac{\delta}2) +\delta =2\delta + 2$, where $(\frac{\delta}2 +
2 + \frac{\delta}2)$ is the length of the path from $x$ to $y$ that
goes via $g$ and $g'$.

So we can find a path of length $\frac{\delta}2 +
(2\delta+2)+\frac{\delta}2=3\delta +2$ from $g$ to $g'$ that lies in
the closed ball of radius $n$.
\end{proof}

It can be shown that if a group is almost convex with constant $C$
with respect to some finite generating set, then the set of all words
that evaluate to the identity and have length at most $C+2$ form a
finite set of relators for a presentation for the group, so the group
is finitely presented. See \cite{\NSnotes}. In this case, however, we 
prove finite presentability directly in the following result.

\begin{thm}\label{thm:subcubicIP}
If the \cg\ for a group $G$ with respect to some finite generating set
$\mathcal G$ has property $L_\delta$ for some $\delta \geq 0$ then $G$ is
finitely presented and has an
isoperimetric function equivalent to  $n^{\frac1{1-\log_3 2}}$.
\end{thm}

\begin{proof}
Assume that the finite generating set $\mathcal G$ is inverse closed.
We will show that every word in $\mathcal G^*$ that evaluates to the
identity has a van Kampen diagram consisting of $2$-cells of perimeter
(at most) $3\delta+2$. Then taking as a set of relators the set of all
words that represent the identity that have length at most
$3\delta+2$, it follows that $G$ is finitely presented.  Moreover, we
will show that the area of a diagram with respect to this presentation
is at most sub-cubic in the length of the word.

 Let $w\in \mathcal G^*$ be a word of length $n$ that evaluates to the
 identity.  We wish to construct a van Kampen diagram for $w$ from
 $2$-cells of perimeter at most $3\delta+2$. So if $n\leq 3\delta +2$
 then it will be a relator, and if $n>3\delta+2$ then proceed as
 follows.

 The word $w$ represents a closed path in the \cg, which starts and
 ends at some vertex $z$. Choose two points $x$ and $y$ that lie at
 distance exactly $\frac{n}{3}$ along the paths $w$ and $w^{-1}$ from
 $z$.  So $x,y$ and $z$ are equally spaced around $w$, that is,
 pairwise at most $ \frac{n}{3}$ apart. Note that $x$ and $y$ need not
 be vertices (they may lie in the interior of edges).  See the left
 side of Figure~\ref{fig:loop}.

Property $L_\delta$ says that there is a point $t$ in the Cayley graph
such that $d(x,t)+d(t,y)\leq d(x,y)+\delta=\frac{n}{3}+\delta$ and
similarly for $x,z$ and $y,z$.

Thus we can find three closed paths containing $x,y,t$, $y,z,t$ and
$x,z,t$ each of perimeter at most
$\frac{n}{3}+(\frac{n}{3}+\delta)=2\frac{n}{3}+\delta$. See the right
side of Figure~\ref{fig:loop}. 

These three closed paths are each strictly shorter than the original
closed path, since $n> 3\delta+2>3\delta$ so $\delta< \frac{n}{3}$ so
$\frac{2n}{3}+\delta < n$.

\begin{figure}
\bt{ccc}
\bt{c}
\includegraphics[width=5.5cm]{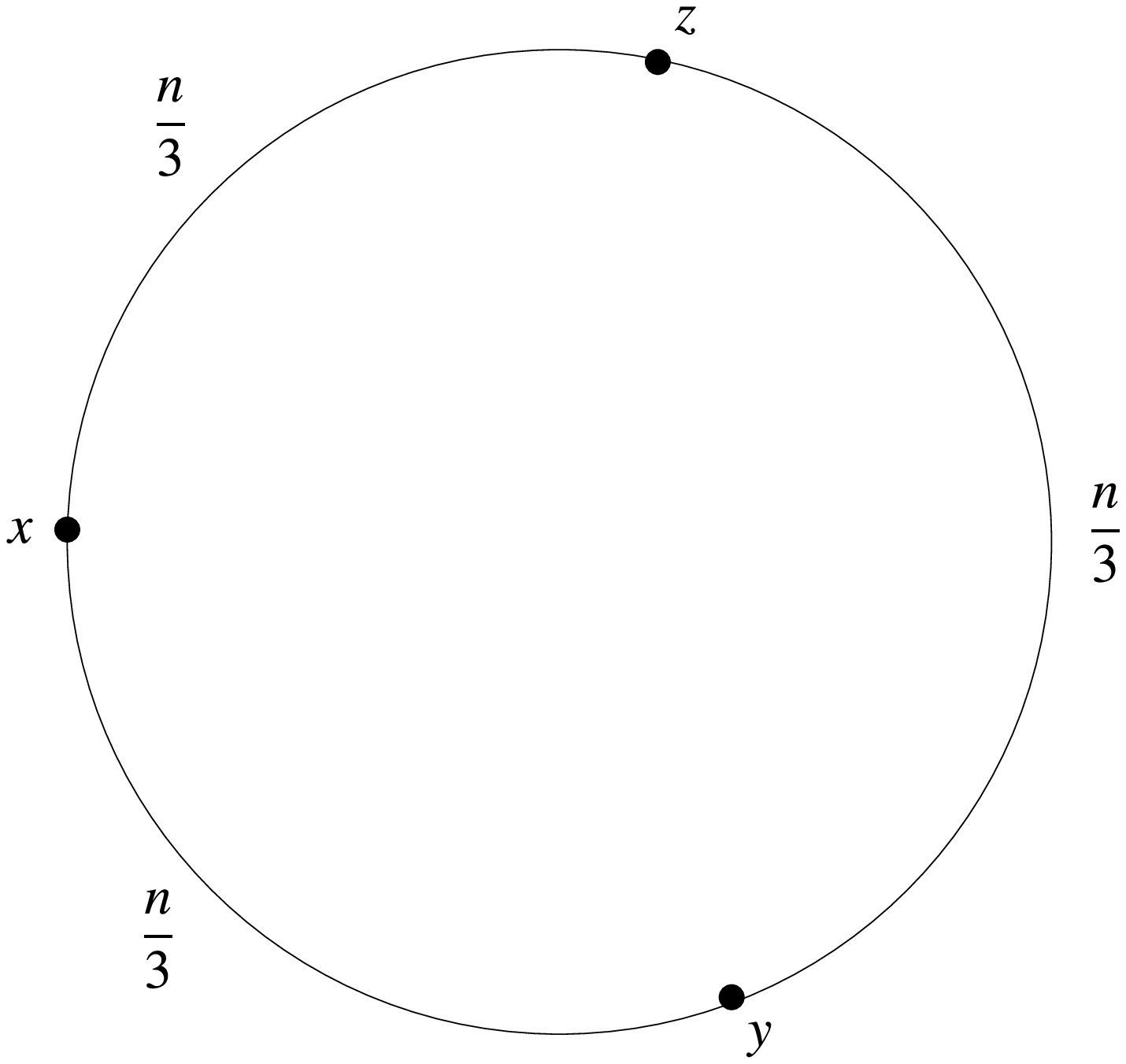}
\et& &
\bt{c}
\includegraphics[width=5.5cm]{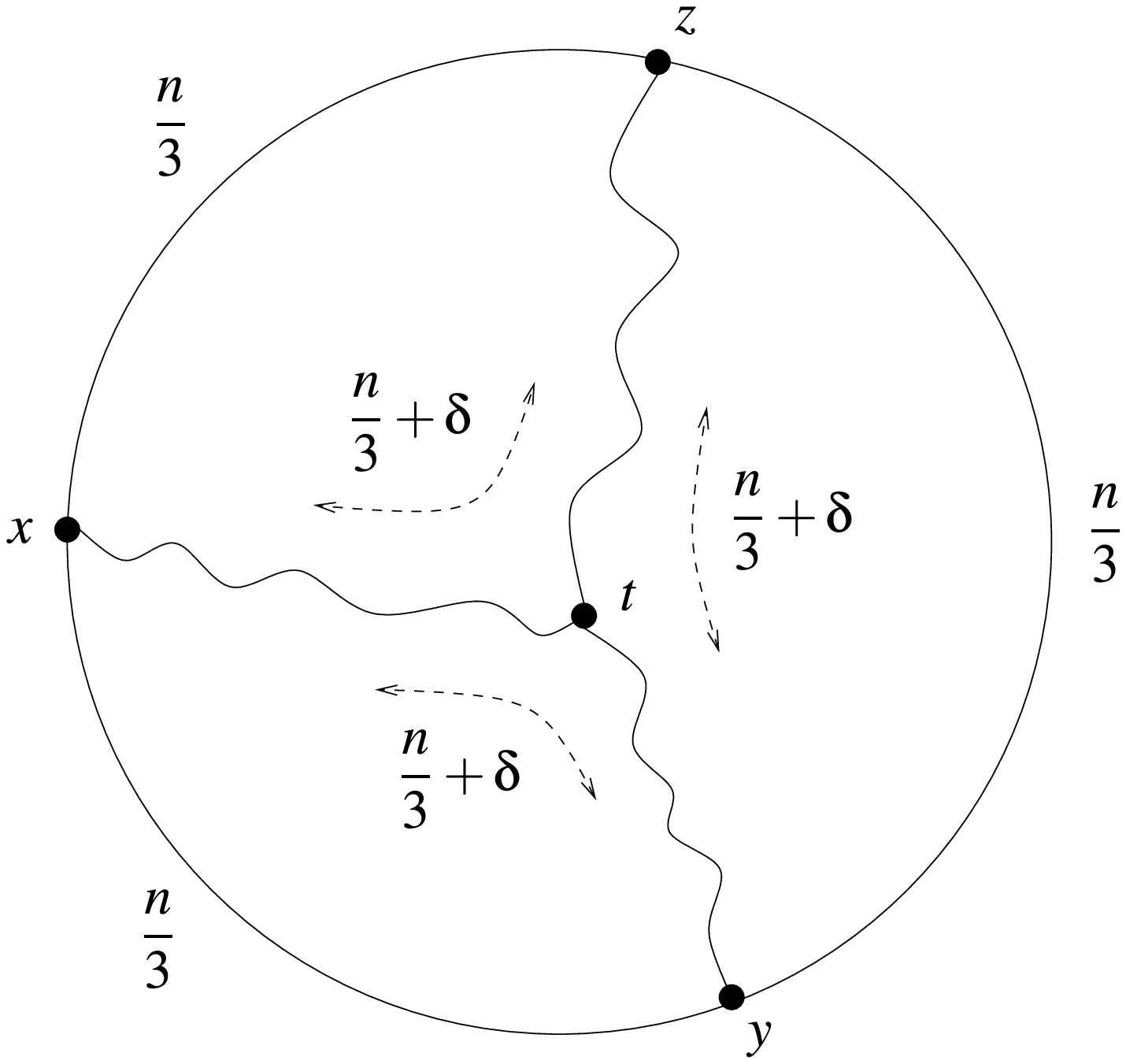}
\et
\et
\caption{Subdividing a closed path into three shorter ones}
\label{fig:loop}
\end{figure}

Note that if $t$ in fact lies on the path $w$ then these paths do not
embed in the \cg.  For the purpose of the argument we are not
concerned with whether these shorter closed paths embed or immerse in the
\cg, we merely want to construct a van Kampen diagram for $w$ with
$2$-cells bounded by closed paths of length at most $3\delta+2$, and
we allow that some of these $2$-cells could be non-embedded.  So we
will iterate this partitioning process until the maximum perimeter of
an internal closed path is not more than $3\delta+2$.

Find $k$ so that $(\frac{3}{2})^k\leq n < (\frac{3}{2})^{k+1}$.  After
one iteration we have three closed paths of length at most
$\frac23n+\delta$.  After a second iteration we get at most nine
closed paths of length at most
$\frac23(\frac23n+\delta)+\delta$
$=(\frac23)^2n+\frac23\delta+\delta$.  After a third iteration
we get at most $3^3$ closed paths of length at most
$\frac23(\frac23(\frac23n+\delta)+\delta)+\delta$
$=(\frac23)^3n+(\frac23)^2\delta+\frac23\delta+\delta$.


Iterating $k$ times we will get at most $3^k$ closed paths of
perimeter at most
\begin{eqnarray*}
&   & (2/3)^kn+(2/3)^{k-1}\delta+\ldots+(2/3)\delta+\delta\\
& < & (2/3)^k(3/2)^{k+1}+\sum_{i=0}^{k-1} (2/3)^i \delta\\
& = &  3/2+\left(\frac{1-(2/3)^k}{1-2/3}\right)\delta\\
& < & 3/2+3\delta.
\end{eqnarray*}

Thus after $k$ iterations we will have partitioned the original closed
path down into closed paths of length less than $3\delta+2$, so we
will have succeeded in finding a van Kampen diagram for $w$ using only
words that evaluate to $1$ and have length less than $3\delta+2$.

Then we have established that $G$ is finitely presented by the
presentation described above.

Now $(\frac{3}{2})^k\leq n$ so taking $\log_3$ of both sides gives
$k\log_3(\frac{3}{2})\leq \log_3 n$.  Let
$c=\frac1{\log_3(3/2)}$. Then $k\leq c\log_3 n$.  So the number of
relators in a van Kampen diagram for $w$ is at most $3^k\leq
3^{c\log_3 n}=(3^{\log_3 n})^c=n^c$.  It follows that the
isoperimetric function is sub-cubic since
$c=\frac1{\log_3(1.5)}=\frac{\ln 3}{\ln 1.5} $ which is approximately
equal to $2.7$. Alternatively, $c=\frac1{\log_3(3/2)}=\frac1{\log_3 3 -
\log_3 2}=\frac1{1 - \log_3 2}$.



Note that the number of iterations required is independent of
$\delta$, so the isoperimetric bound is not improved by smaller
$\delta$ (such as $0$).
\end{proof}

\section{Remarks}

The 3-dimensional integral Heisenberg group has a cubic isoperimetric function (see
\cite{\Epstein} p.165), so in the light of Theorem
\ref{thm:subcubicIP} cannot enjoy property $L_\delta$ for any~$\delta\geq 0$.


An open question is whether there is a group with an isoperimetric
function greater than quadratic and less than the sub-cubic bound
given in Theorem \ref{thm:subcubicIP} that has $L_{\delta}$. Brady and
Bridson have a family of groups having isoperimetric functions with
exponent $2\log_2 \frac{2p}{q}$ for all $p\geq q$ \cite{\BB}, however a quick
investigation of these suggests that they would not enjoy $L_{\delta}$.


\def\cprime{$'$}


\Addresses
\recd
\end{document}